\documentclass[12pt]{article}
\usepackage{amsmath}
\usepackage{amssymb}
\usepackage{amsfonts}
\usepackage{eucal}
\usepackage[usenames]{color}
\usepackage{graphicx}
\numberwithin{equation}{section}
\newtheorem{Theorem}{Theorem}[section]
\newtheorem{Proposition}[Theorem]{Proposition}

\newtheorem{Remark}[Theorem]{Remark}

\begin{document}
\title{Some examples of solutions to
an inverse problem for the first-passage place of a jump-diffusion process.
}

\author{Mario Abundo\thanks{Dipartimento di Matematica, Universit\`a  ``Tor Vergata'', via della Ricerca Scientifica, I-00133 Rome, Italy.
E-mail: \tt{abundo@mat.uniroma2.it}}
}
\date{}
\maketitle

\begin{abstract}
\noindent We report some additional examples of explicit solutions to an inverse first-passage place problem for one-dimensional diffusions with jumps, introduced in a previous paper. If $X(t)$ is a one-dimensional diffusion with jumps, starting
from a random position $\eta \in [a,b],$ let be $\tau_{a,b}$ the time at which $X(t)$ first exits the interval $(a,b),$
and $\pi _a = P(X(\tau_{a,b}) \le  a)$ the probability of exit from the left of $(a,b).$
Given a probability $q \in (0,1),$
the  problem consists in finding the density $g$ of $\eta$ (if it
exists) such that $\pi _a = q;$ it can be seen as a problem of optimization.
\noindent

\end{abstract}


\section{Introduction and preliminary results}
This short note is a continuation of the paper \cite{abundo:saa20}, in which we have studied an inverse first-passage place problem (IFPP) for a one-dimensional jump-diffusion process; while for simple-diffusions (i.e. without jumps) a number of examples was reported in \cite{abundo:saa20}, in that article we were able to present only one example, concerning diffusions with jumps. Therefore, in this paper we report
some additional  examples of explicit solutions to the IFPP problem for diffusions with jumps.
\par
The IFPP problem, as well as the analogous inverse first-passage time problem, have interesting applications in Mathematical Finance, in particular in credit risk modeling, where
the first-passage time represents a default event of an obligor (see e.g.
\cite{jackson:stapro09}), in Biology,
in the scope of diffusion models for neural
activity (see e.g. \cite{lanska:89}), Engineering, and  many other fields; for more about inverse first-passage time problems, see e.g
\cite{abundo:saa19}, \cite{abundo:mathematics},  \cite{abundo:LNSIM}, \cite{abundo14b}
\cite{abundo13b}, \cite{abu:smj13}, \cite{abundo:stapro13}, \cite{abundo:stapro12}, \cite{jackson:stapro09}.
As regards the direct first-passage time problem for jump-diffusions, see e.g. \cite{abundo:mcap10}, \cite{abundo:pms00}, \cite{kouwang:03}, \cite{tuckwell:76}); as for the direct first-passage place problem, few results are known:
it was studied by Lefebvre (\cite{lefeb:20}, \cite{lefeb:19}, \cite{lefeb:19b}), and Kou and Wang (\cite{kouwang:03}), where
equations for the moments of first-passage places were established; in a particular case Lefebvre
found exact formulae for the $k$-th moments of the first-passage place, providing also an approximate analytical expression for it.
 \par
We recall the terms of the IFPP problem.
Let be
\begin{equation} \label{jumpdiffueq2}
X(t) =  \eta + \int _0 ^t \mu(X(s))ds + \int _0 ^t \sigma(X(s))dB_s + \sum _{i=1} ^{N_1 (t)} \varepsilon _i (X(t)) +  \sum _{i=1} ^{N_2 (t)} \Delta _i (X(t)),
\end{equation}
a one-dimensional, time-homogeneous jump-diffusion
process starting from a random position $\eta \in [a,b] ,$
where $B_t$ is standard Brownian motion,  $\mu (\cdot)$ and $ \sigma ( \cdot)$ are smooth enough deterministic functions, and $\{N_k (t) \}$ is a time-homogeneous  Poisson process with rate $\lambda_k >0,$ for $k=1, 2.$
The three stochastic processes $B_t, \ N_1 (t)$ and $N_2 (t)$ are assumed to be independent, and the r.v. $\eta$ is independent of them; moreover the state-dependent
random variables $\varepsilon _i(X(t)) >0$ and $\Delta _i(X(t)) <0, \ i=1,2, \dots ,$ are independent and identically distributed, and are independent between themselves. Note that Eq. \eqref{jumpdiffueq2} is slightly different
from  the analogous equation in \cite{abundo:saa20}, which defines the jump-diffusion process $X(t),$ there considered; really,  the representation of $X(t)$ by means of \eqref{jumpdiffueq2} allows to points out explicitly the
positive and negative jumps; this formulation was inspired by Lefebvre's paper \cite{lefeb:20}. \par\noindent
We suppose that the first-exit time of $X(t)$ from the interval $(a,b),$ namely
\begin{equation}
\tau _ {a, b} = \inf \{t \ge 0: X(t) \notin (a, b) \},
\end{equation}
is finite with probability one,
and let $X( \tau _{a,b})$ be the first-passage place of $X(t)$ at time $\tau _{a,b}.$ By assumption, one has
$X( \tau _{a, b}) \le a $ or $X(\tau _{a, b}) \ge b;$ we denote by  $\pi _a = P ( X(\tau _{a, b}) \le a )$ the probability that the process $X(t)$ first exits the interval $(a,b)$
from the left, and
by $\pi _b = 1- \pi _a = P ( X(\tau _{a, b} ) \ge b )$ the probability  that $X(t)$  first exits from the right. \par\noindent
Actually,
we have considered in  \cite{abundo:saa20} the following inverse
first-passage place (IFPP) problem:
\bigskip

{\it given a probability $q \in (0,1),$ find the density $g$ of $\eta$ (if it
exists) for which it results \par
%$P(Z = a)= q.$
$\pi_a = q .$ }
\bigskip

\noindent The function $g$ is called a solution to the IFPP
problem. In fact, the solution to the IFPP problem, if it exists, is not necessarily unique (see \cite{abundo:saa20}).
As we will see in Remark \ref{remopt}, the IFPP problem can be also seen as a problem of optimization.
\bigskip

In the next Section, we present some additional examples of explicit solutions to the IFPP problem for one-dimensional diffusions $X(t)$ with jumps.
They also provide information about the corresponding direct first-passage time problem, since they involve the calculation of the exit probability of $X(t)$
from the left  of the interval $(a,b).$
\par
Let $f_\varepsilon (\epsilon)$ and $f_\Delta (\delta)$ be the probability density functions of the random variables $\varepsilon _i ( \cdot ) >0$ and $\Delta _i ( \cdot ) <0,$ respectively; we suppose that
the infinitesimal drift $\mu(x)$ and the infinitesimal diffusion coefficient  $\sigma (x)$ of the process $X(t)$ are smooth enough deterministic functions, and we denote by
$\tau _{a, b}(x)$ the first-exit time of $X(t)$ from the interval $(a,b),$ with the condition that  $\eta =x \in [a,b].$ Moreover, we set
$\pi _a (x) = P(X(\tau _{a ,b} (x)) \le a)$ and $ \pi_b (x) = P(X(\tau _{a, b}(x)) \ge b)= 1 - \pi _b (x).$ \par
We recall (see e.g. \cite{abundo:saa20}, \cite{lefeb:20}, \cite{lefeb:19}) that
the function $v(x) := \pi _a (x)$
satisfies the integro-differential problem with outer conditions:
\begin{equation} \label{eqforpia}
\begin{cases}
\frac 1 2 \sigma ^2 (x) v'' (x) + \mu (x) v' (x) + \lambda _1 \int _ {- \infty } ^ {+ \infty} [v(x+ \epsilon) - v(x)] f_ \varepsilon (\epsilon ) d \epsilon + \\
\ + \lambda _2 \int _ {- \infty } ^ {+ \infty} [v(x+ \delta) - v(x)] f_ \Delta (\delta ) d \delta
 = 0, \ x \in (a,b) \\
v(x)= 1 \ {\rm if} \ x \le a \ {\rm and } \ v(x)= 0 \ {\rm if } \ x \ge b.
\end{cases}
\end{equation}
If there is no jump, that is, $f_\varepsilon (\epsilon)$ and $f_\Delta (\delta)$ are identically zero, then the process defined by
\eqref{jumpdiffueq2} is a (continuous) simple-diffusion, and so the outer conditions in \eqref{eqforpia} become the boundary conditions
$v(a)=1, \ v(b)=0.$
\bigskip

Returning back to the case when the jump-diffusion  $X(t)$ starts from the random position $\eta \in [a,b],$ we
suppose that $\eta$ possesses a density $g(x);$ then the following holds (see \cite{abundo:saa20}):
\begin{Proposition} \label{mainpro}
Let $X(t)$ be the jump-diffusion process defined by \eqref{jumpdiffueq2};
with the previous notations,
if a solution $g$ exists to the IFPP problem for $X(t)$ and $q \in (0,1),$ then the function $g$ must satisfy the following equation:
\begin{equation} \label{maineq}
q= \int _a ^b g(x) \pi_a (x) dx ,
\end{equation}
where $\pi_a (x)$ is the solution of \eqref{eqforpia}.
\end{Proposition}
\hfill  $\Box$
\bigskip
\begin{Remark}
For an assigned $q \in (0,1),$ Eq. \eqref{maineq} is an integral equation in the
unknown $g(x).$ Unfortunately, no method is available to solve analytically this equation, so
any possible solution $g$ to the IFPP problem  must be found by making attempts (see also Remark 2.5 in \cite{abundo:saa20}).
\end{Remark}

\begin{Remark} \label{remopt}
The IFPP problem can be seen as a problem of optimization: indeed, let  ${\cal G}$ be the set of  probability densities on the interval
$(a,b),$ and consider the functional $\Psi : {\cal G} \longrightarrow \mathbb{R} ^ +$ defined, for any $g \in {\cal G}$, by
\begin{equation}
\Psi (g) = \left ( q- \int _ a ^b g(x) \pi_a(x) dx \right ) ^2 .
\end{equation}
Then, a solution $g$ to the IFPP problem, is characterized by
\begin{equation}
g = \arg \min _ { g \in {\cal G} } \Psi (g) .
\end{equation}
Of course, if there exists more than one function $g \in {\cal G}$ at which $\Psi (g)$ attains the minimum, the solution of the IFPP problem is not unique.
\end{Remark}
\bigskip

\section{Examples}
\noindent {\bf Example 1.} Let $X(t)$ be a jump-diffusion process of the form \eqref{jumpdiffueq2}; we suppose that $\varepsilon _i (X(t)),$ given that $X(t) = \xi,$ is uniformly distributed on the interval $(0, \alpha _1 \xi ), $ where $ \alpha _1 >0;$
in analogous way, we assume that $\delta_i (X(t)),$ given that $X(t) = \xi,$ is uniformly distributed on the interval $(-\alpha _2 \xi, 0 ), $
with $ 0 < \alpha _2  \le 1 .$
Moreover, we suppose that the drift is $\mu (x) = \frac 1 2 (\lambda _2 \alpha _2 - \lambda _1 \alpha _1 )x ,$
while $\sigma (x)$ is any diffusion coefficient,
and let be $\alpha, \ \beta $ positive constants; then, a solution
 $g$ to the IFPP problem for $X(t)$ and  $q= \frac \beta { \alpha + \beta }$ is the modified Beta density in the interval $(a,b)$ with parameters $\alpha$ and $\beta ,$ namely:
\begin{equation} \label{solution1}
g(x) = \frac 1 {(b-a)^{ \alpha + \beta -1} } \cdot \frac {(x-a) ^{ \alpha -1} (b-x) ^{ \beta -1} } {B( \alpha, \beta) } \cdot \mathbb {I} _ {(a,b)}(x),
\end{equation}
where $B( \alpha, \beta ) = \frac { \Gamma ( \alpha ) \Gamma ( \beta )} { \Gamma ( \alpha + \beta )}$
(the ordinary Beta density is obtained for $a=0$ and $b=1).$ \par\noindent
In fact, from \eqref{eqforpia} the equation for $\pi _a (x)$ is
 (see also \cite{lefeb:20}):
 \begin{equation} \label{eqpialef}
 \frac 1 2 \sigma ^2 (x) v''(x) + \mu (x) v'(x) - (\lambda _1 + \lambda _2) v(x) + \frac {\lambda _1 } { \alpha _1 x} \int_0^ {\alpha _1 x} v(x+ \epsilon) d \epsilon + \frac {\lambda _2 } { \alpha _2 x} \int_{- \alpha _2 x} ^0 v(x+ \delta) d \delta =0,
 \end{equation}
with the conditions
\begin{equation}
 v(x)= 1, \ x \le a; \ v(x) = 0, \ x \ge b ,
 \end{equation}
and it is satisfied by
\begin{equation}
\pi _a (x)=
\begin{cases}
1 & \ {\rm if} \ x \le a \\
\frac {b-x } {b-a } & \ {\rm if } \ x \in (a,b) \\
0 & \ {\rm if} \ x \ge b ,
\end{cases}
\end{equation}
irrespective of the diffusion coefficient $ \sigma (x).$
Then, to verify that $g,$ given by \eqref{solution1}, is solution to the IFPP problem, it suffices to substitute $g, \ q$ and $\pi _a (x)$ into Eq. \eqref{maineq}
(in the calculation of the integral, one can use that the mean of the r.v. $\eta$ with density $g$ is $(a \beta + b \alpha)/( \alpha + \beta );$ in fact,
one has $\eta = a + (b-a) U,$ being  $U$ a r.v. with Beta density).
 \par\noindent
Note that, for $\beta > \alpha$ it results $q > 1/2,$ if $ \beta = \alpha $ one has $ q=1/2,$ while for $\beta < \alpha $ one has $ q < 1/2.$
For $\alpha = \beta =1, \ g$ turns out to be the uniform density in the interval $(a,b).$  \par
We remark that the simple-diffusion process $\widetilde X(t)$ obtained by $X(t)$ disregarding the jumps (that is, setting $f_ \varepsilon (\epsilon) = f_ \Delta ( \delta ) =0),$
is driven by the SDE
\begin{equation}
d\widetilde X(t)= \frac 1 2 (\lambda _2 \alpha _2 - \lambda _1 \alpha _1 ) \widetilde X(t) dt + \sigma (\widetilde X(t)) dB_t .
\end{equation}
 Thus,
$\widetilde X(t)$ is (see also \cite{lefeb:20}): \bigskip

\noindent
$\bullet $ Brownian motion, if $\lambda _2 \alpha _2 = \lambda _1 \alpha _1$ and $\sigma (x) =1,$  \par\noindent
$ \bullet $ Ornstein-Uhlenbeck process, if
$\lambda _2 \alpha _2 < \lambda _1 \alpha _1$ and $\sigma (x) = const.,$ \par\noindent
$ \bullet $ Geometric Brownian motion, if
$\lambda _2 \alpha _2 > \lambda _1 \alpha _1 $ and $\sigma (x) = c x,$ with $c$ a positive constant, \par\noindent
$ \bullet $ the CIR-like model in mathematical finance, if $\sigma (x) = \sqrt { x \vee 0} ,$ \par\noindent
$\bullet $ the Wright$\&$Fisher-like process, if $\sigma (x) = \sqrt { x(1-x) \vee 0} $ (see e.g. \cite{abundo:saa20}). \bigskip

Note that
$\frac {b-x } {b-a }, \ x \in (a,b),$ is nothing but the exit probability of $X(t)= x + B_t$ at the left of the interval $(a,b);$  thus,
the function $g$ given by \eqref{solution1} is also solution to the IFPP problem for Brownian motion and $q= \frac \beta { \alpha + \beta }, \ \alpha , \ \beta >0$ (see Example 3 of
\cite{abundo:saa20}).
\bigskip

\noindent
{\bf Example 2.} Take $a=0, \ b=1, \ \gamma >0,$ and suppose that $X(t)$ is the jump-diffusion \eqref{jumpdiffueq2} with
diffusion coefficient $ \sigma (x) = \sqrt {x \vee 0} $ and
linear drift
$\mu (x) = Ax + B,$
where
$$ A= \frac 1 \gamma  \ \left [ \lambda _1 + \lambda _2 + \frac 1 { \gamma +1 } \left ( \frac { \lambda _1} {\alpha _1} (1- (1 + \alpha _1) ^{\gamma +1} )
+  \frac { \lambda _2} {\alpha _2} ((1 - \alpha _2) ^{\gamma +1} -1) \right )   \right ], \ B= - \frac 1 2 ( \gamma -1).$$
We assume that the functions $f_ \varepsilon (\epsilon), \ f_ \Delta ( \delta ) $ are the same ones, as in Example 1. Then, for positive $\alpha, \ \beta,$ a solution $g$ to the
IFPP problem for $X(t)$ and $q= 1 - \frac { \Gamma ( \alpha + \gamma)\Gamma ( \alpha + \beta)} {\Gamma (\alpha) \Gamma ( \alpha + \beta + \gamma) } $ is the Beta density in $(0,1)$ with parameters $\alpha, \beta,$  that is
$g(x)= \frac {\Gamma ( \alpha + \beta ) } {\Gamma ( \alpha) \Gamma (\beta) } \cdot x ^{ \alpha -1} (1-x) ^{ \beta -1} \cdot \mathbb {I} _ {(0,1)}(x) .$ \par\noindent
In fact, for the above infinitesimal coefficients $\mu(x)$ and $\sigma (x),$ it is easy to see that Eq. \eqref{eqpialef} is satisfied by  $\pi _0 (x) = 1-x^ \gamma$ for $x \in (0,1),$ and so to verify
that $g$ is solution to the IFPP problem, it is enough
to substitute $g, \ q$ and $\pi _0 (x)$ into Eq. \eqref{maineq} (now, the integral in \eqref{maineq} is nothing but $1-E(Z^ \gamma ),$ where $Z$ is a r.v.
with Beta density; thus, in the calculation it is convenient to use that
$E(Z ^ \gamma )= \frac {\Gamma (\alpha + \beta ) \Gamma (\alpha + \gamma )} {\Gamma (\alpha ) \Gamma (\alpha + \beta + \gamma)}$ (see e.g. \cite{gupta})).
\par\noindent
For instance, if $\gamma =2,$ then $\pi _0 (x) = 1-x^2,$ and $q= \frac { \beta (\beta + 2 \alpha +1)} { (\alpha + \beta )(\alpha + \beta +1)} .$  \par\noindent
Notice that the simple-diffusion process $\widetilde X(t),$ obtained by $X(t)$ disregarding the jumps, satisfies the SDE:
\begin{equation}
d \widetilde X(t)= (A \widetilde X(t) +B) dt + \sqrt { \widetilde X(t) \vee 0} \ dB_t ,
\end{equation}
 which provides a special case of
the CIR model. \par\noindent
Really, $\pi_0(x)= 1-x^2$ is also the exit probability at the left of $(0,1)$ of the simple-diffusion driven by the SDE:
\begin{equation} \label{SDEa}
d \widetilde X(t)= - \frac 1 2 dt + \sqrt {\widetilde X(t) \vee 0} \ dB_t .
\end{equation}
Thus,
the Beta density in $(0,1)$ is also solution to the IFPP problem for the diffusion driven by \eqref{SDEa} and $q= \frac { \beta (\beta + 2 \alpha +1)} { (\alpha + \beta )(\alpha + \beta +1)} .$
\bigskip

\noindent {\bf Example 3.} Take $a=0, \ b=1, \ \gamma >0$ and suppose that $X(t)$ is the jump-diffusion \eqref{jumpdiffueq2} with $\sigma (x) = \sqrt {x(1-x) \vee 0} ,$
$\mu (x) = A' x + B  ,$  with $A' = \frac 1 2 ( \gamma -1) + A,$ and the constant $A, \ B,$ as well as
the functions $f_ \varepsilon (\epsilon), \ f_ \Delta ( \delta ) ,$ are the same ones as in Example 2. Then, for positive $\alpha, \ \beta,$ a solution $g$ to the
IFPP problem for $X(t)$ and $q= 1 - \frac { \Gamma ( \alpha + \gamma)\Gamma ( \alpha + \beta)} {\Gamma (\alpha) \Gamma ( \alpha + \beta + \gamma) } $ is the Beta density in $(0,1)$ with parameters $\alpha, \beta .$
\par\noindent
In fact, for the above infinitesimal coefficients  Eq. \eqref{eqpialef} is satisfied by $\pi _0 (x) = 1-x^\gamma$ for $x \in (0,1),$ as in Example 2; thus to verify the result it is enough
to substitute $g, \ q$ and $\pi _0 (x)$ into Eq. \eqref{maineq}.
\par\noindent
Notice that the simple-diffusion process $ \widetilde X(t),$ obtained by $X(t)$ disregarding the jumps, satisfies the SDE:
\begin{equation}
d \widetilde X(t)= (A' \widetilde X(t) + B) dt + \sqrt { \widetilde X(t)(1- \widetilde X(t)) \vee 0} \ dB_t ,
\end{equation}
 which provides the Wright-Fisher-like process
(see e.g. \cite{abundo:saa20}).
\par\noindent
Really, $\pi _0 (x) = 1- x^2$
is also the exit probability at the left of $(0,1)$ of the simple-diffusion driven by the SDE:
\begin{equation} \label{SDEb}
d \widetilde X(t)= \frac 1 2 (\widetilde X(t)-1) dt + \sqrt { \widetilde X(t)(1- \widetilde X(t)\vee 0 } \ dB_t .
\end{equation}
Thus,
the Beta density in $(0,1)$ is also solution to the IFPP problem for the simple-diffusion driven by \eqref{SDEb} and $q= \frac { \beta (\beta + 2 \alpha +1)} { (\alpha + \beta )(\alpha + \beta +1)}.$
\bigskip

\noindent {\bf Example 4.} With the previous notations and assumptions on the Poisson processes $N_k(t),$ let be
$ \bar \epsilon, \ \bar \delta >0,$ and suppose that, for $\eta \in [a,b]:$
\begin{equation}
X(t) = \eta + (\bar \delta \lambda _2 - \bar \epsilon \lambda _1)t + \int _0 ^t \sigma (X(s)) dB_s + \bar \epsilon  N_1(t) - \bar \delta  N_2 (t) .
\end{equation}
Then, a solution
 $g$ to the IFPP problem for $X(t)$ and  $q= \frac \beta { \alpha + \beta } \ (\alpha, \ \beta >0),$ is the modified Beta density in the interval $(a,b),$
 given by \eqref{solution1}.
 In fact, now the equation for $v(x)= \pi _a (x)$ becomes:
 \begin{equation}
 \frac 1 2 \sigma ^2(x) v'' (x) + (\bar \delta \lambda _2 - \bar \epsilon \lambda _1) v'(x) - (\lambda _1 + \lambda _2 ) v(x) + \lambda _1 v(x+ \bar \epsilon) + \lambda _2 v(x - \bar \delta) =0,
 \ x\in (a,b),
 \end{equation}
which is satisfied by $v(x)= \pi _a (x) = \frac {b-x } {b-a }, \ x \in (a,b),$ irrespective of $\sigma (x);$  thus, the assertion is soon verified, proceeding as in Example 1. \par
A variant is obtained by considering the jump-diffusion:
\begin{equation}
X(t)= \eta - \bar \epsilon \lambda _1 t + \int _0 ^t \sigma (X(s)) dB_s + N_1(t);
\end{equation}
then,  a solution
 $g$ to the IFPP problem for $X(t)$ and  $q= \frac \beta { \alpha + \beta } \  (\alpha, \ \beta >0),$ is again the modified Beta density in the interval $(a,b).$
 It suffices to note that now
 the equation for $v(x)= \pi _a (x)$ is
 \begin{equation}
  \frac 1 2 \sigma ^2(x) v'' (x) - \bar \epsilon \lambda _1 v'(x) - \lambda _1 v(x) + \lambda _1 v(x+ \bar \epsilon) =0, \ x\in (a,b),
 \end{equation}
and it is satisfied again by  $\pi _a (x) = \frac {b-x } {b-a }, \ x \in (a,b),$ irrespective of $\sigma (x).$
\bigskip

\noindent {\bf Example 5.} Take $a=0, \ b=1;$ for $ \bar \epsilon, \ \bar \delta >0,$
suppose that:
\begin{equation}
dX(t)= \mu (X(t))dt + \sqrt { X(t) } dB_t + \bar \epsilon dN_1(t) - \bar \delta dN_2 (t), \ X(0) = \eta \in [0,1] ,
\end{equation}
where
\begin{equation} \label {muex5}
\mu (x)= \frac 1 { \ln 2 } \left [ - \frac {(\ln 2)^2 x} 2 - \lambda _1 (2^ { \bar \epsilon } -1 ) + \lambda _2 (1- 2 ^ {- \bar \delta }) \right ] ,
\end{equation}
and $N_k(t)$ are Poisson processes with intensity $\lambda _k, \ k=1,2$ (we can write $\sqrt {X(t)}$ instead of
$\sqrt {X(t)  \vee 0},$ since $X(t)$ is  $\ge 0$
until  the first-exit time of $X(t)$ from the interval
$(0,1)).$
Then, for any $\alpha, \beta >0$ a solution
 $g$ to the IFPP problem for $X(t)$ and
 \begin{equation}
 q= 2 - \sum _{k=0} ^ \infty \frac {(\ln 2 )^k } { k!}  \ \frac {B(\alpha +k, \beta ) } {B(\alpha, \beta ) } ,
 \end{equation}
 is the Beta density in the interval $(0,1)$ with parameters $\alpha$ and $\beta $
(if e.g. $\alpha = \beta =1,$ one has $q= 2 - 1/ \ln 2$ and $g$ is the uniform density in $(0,1),$ if
$\alpha = \beta =2,$ then $q= 2 - 6 \cdot \frac {3 \ln 2 -2 } {(\ln 2 )^3 } $ and $g(x)= 6x(1-x), \ x\in (0,1) ).$ \par\noindent
 In fact, now the equation for $v(x)= \pi _0 (x)$ becomes:
 \begin{equation}
 \frac 1 2 x v'' (x) + \mu (x) v'(x) - (\lambda _1 + \lambda _2 ) v(x) + \lambda _1 v(x+ \bar \epsilon) + \lambda _2 v(x - \bar \delta) =0,
 \ x\in (0,1),
 \end{equation}
which is satisfied by $v(x)= \pi _ 0 (x) = 2 - 2^x , \ x \in (0,1),$ if $\mu (x)$ is given by \eqref{muex5};
thus, to verify
that $g (x)= \frac {\Gamma ( \alpha + \beta ) } {\Gamma ( \alpha) \Gamma (\beta) }  x ^{ \alpha -1} (1-x) ^{ \beta -1} $ is solution to the IFPP problem, it is enough to substitute $g, \ q$ and $\pi _0 (x)$ into Eq. \eqref{maineq}
(to calculate the integral in \eqref{maineq} it is convenient to note that
\begin{equation}
\int _0 ^1 2^x g(x) \ dx = E \left ( e^ {(\ln 2) X} \right ),
\end{equation}
where
\begin{equation}
E \left ( e^ {tX} \right ) = \sum _{k=0} ^ \infty \frac { t^k } { k!}  \ \frac {B(\alpha +k, \beta ) } {B(\alpha, \beta ) }
\end{equation}
is the moment generating function of a r.v. having Beta density, with parameters $\alpha$ and $\beta ,$ see e.g \cite{gupta}).
\bigskip

\noindent {\bf Example 6.}  Take $a=0, \ b =1$ and, for $ \eta \in [0,1]$ consider the jump-diffusion:
\begin{equation}
X(t) = \eta - \frac \pi 4 \int _0 ^t \cos \left ( \frac \pi 2 X(s) \right ) ds + \int _0 ^t \sqrt { \sin \left ( \frac \pi 2 X(s) \right )} \ dB_s + 4 N_1(t),
\end{equation}
where $N_1(t)$ is a Poisson Process with intensity $\lambda _1$ (note that the amplitude of jumps is $4$ and, as soon as a jump occurs, the process
 exits $(0,1)$ from the right). \par\noindent
Then,  a solution
 $g$ to the IFPP problem for $X(t)$ and  $q= 2 / \pi ,$ is the uniform density in the interval $(0,1).$
Now, the equation for $v(x)= \pi _0 (x)$ becomes:
 \begin{equation}
  \frac 1 2 \sin \left ( \frac \pi 2 x \right ) v'' (x) - \frac \pi 4 \cos \left ( \frac \pi 2 x \right ) v'(x) - \lambda _1 v(x) + \lambda _1 v(x+ 4) =0, \ x\in (0,1),
 \end{equation}
which is satisfied by $v(x)= \pi _0 (x)= \cos \left ( \frac \pi 2 x  \right ),$ for $x \in (0,1).$ Thus, to verify the result it is enough
to substitute $g (x)= \mathbb{I}_{(0,1)}(x),  \ q$ and $\pi _0 (x)$ into Eq. \eqref{maineq}. \par\noindent
Note that $\pi _0 (x)= \cos \left ( \frac \pi 2 x  \right )$ is also the exit probability at the left of $(0,1)$ of the simple-diffusion obtained from $X(t)$ disregarding the jumps.
\bigskip

Finally, we recall the following example, already presented in \cite{abundo:saa20}, in which $a=0, \ b= 2 \epsilon \ (\epsilon $ a fixed positive number), and the exit probability, $\pi _0 (x),$ from the left of the interval $(0, b)$ has a more complicated form, since it is not a polynomial, exponential-like, or  trigonometric function.\bigskip

\noindent {\bf Example 7.}
For $\epsilon >0,$ take $a=0, \ b= 2 \epsilon ,$ and let be $X(t)= \eta + B_t + \epsilon N_1(t),$
where the
starting point $\eta $ is random in $[0,2 \epsilon]$ and
$N_1(t)$ is a time-homogeneous Poisson process with rate $\lambda_1 =1.$
Conditionally to $\eta =x \in [0, 2 \epsilon],$  the probability
$\pi _0 (x) = P(X( \tau _{0 , 2 \epsilon}(x) ) \le 0 ) $
is the solution to the integro-differential problem:
\begin{equation} \label{eqforpizero}
\begin{cases}
\frac 1 2  v'' (x)  +  v(x+\epsilon) - v(x), \ x \in (0, 2 \epsilon) \\
v(x)= 1 \ {\rm if} \ x \le 0 \ {\rm and } \ v(x)= 0 \ {\rm if } \ x \ge 2 \epsilon.
\end{cases}
\end{equation}
The computation of  $\pi _0 (x)$ is very complicated, however its explicit form was found in \cite{abundo:pms00}, and  it is given by:
\begin{equation} \label{pi0forexampleuno}
\pi _0 (x) =
\begin{cases}
e^{- x \sqrt 2 }(A + ax) + e ^{x \sqrt 2 } (B+ bx), \ x \in (0, \epsilon) \\
c \left ( e^{ - x \sqrt 2 } - e ^{-4 \epsilon \sqrt 2 + x \sqrt 2 } \right ), \ x \in [\epsilon, 2 \epsilon)
\end{cases}
\end{equation}
where $a, \ b, \ c, \ A$ and $B$ are constants such that
\begin{equation}
 a = \frac {ce^{- \epsilon \sqrt 2} } { \sqrt 2}, \ b= - \frac {c e^{ -3 \epsilon \sqrt 2 }} { \sqrt 2 }, \ B=1-A,
 \end{equation}
and $A, \ c$ have to be found by requiring that  $\pi_0 (x)$ is a $C^2$ function; doing this, one finds:
\begin{equation}
A = \frac {e^{\epsilon \sqrt 2}(\beta \sqrt 2 - \delta) } {\alpha \delta - \beta \gamma } , \ c = \frac {e^{  \epsilon \sqrt 2} (\gamma - \alpha \sqrt 2) } {\alpha \delta - \beta \gamma },
\end{equation}
where
\begin{equation}
\alpha = -2 \sinh ( \epsilon \sqrt 2 ), \ \beta= e^{- \epsilon \sqrt 2} (e^{-2 \epsilon \sqrt2 } -1) ,
\end{equation}
\begin{equation}
\gamma = -2 \sqrt 2 \cosh ( \epsilon \sqrt 2 ), \ \delta = - 2 \epsilon e^{ -2 \epsilon \sqrt 2 } + \sqrt 2 (e ^{ - \epsilon \sqrt 2 } + e ^{-3 \epsilon \sqrt 2 } ) .
\end{equation}
Thus, $\pi_0 = \int _0 ^{2 \epsilon} \pi _0 (x) g(x) dx $ can be  calculated, after tedious calculations, for any explicit density $g(x), \ x \in (0, 2 \epsilon).$ \par
Now, let be
$$ q = \frac 1 { 2 \epsilon \sqrt 2} \left [ \frac {\gamma - (\alpha + \beta ) \sqrt 2 + \delta } {\alpha \delta - \beta \gamma }
+ e^{\epsilon \sqrt 2 } \left (1- \frac {e^{\epsilon \sqrt 2} (\beta \sqrt 2 - \delta )  } { \alpha \delta - \beta \gamma} \right ) \right ] $$
$$ +  \frac 1 { 2 \epsilon \sqrt 2 } \left [ \frac {\sqrt 2 e^{- \epsilon \sqrt 2} \left (\sqrt 2 - \epsilon - \frac 1 4 \right ) ( \gamma -  \alpha \sqrt 2 ) } {\alpha \delta - \beta \gamma }
 \right ]$$
\begin{equation} \label{qofexample1}
+ \frac 1 { 2 \epsilon \sqrt 2 } \left [  \frac {e^{-2 \epsilon \sqrt 2} ( \gamma - \alpha \sqrt 2 ) } {\alpha \delta - \beta \gamma }
+ \frac{2 e^ {\epsilon \sqrt 2 } ( \beta \sqrt 2 - \delta ) } {\alpha \delta - \beta \gamma  }   +  \frac {\gamma - \alpha \sqrt 2  } {2(\alpha \delta - \beta \gamma  )} -1 \right].
\end{equation}
Then, a solution $g$ to the IFPP problem for $X(t)$ and the above value of $q$ is the uniform density in $(0,2 \epsilon )$ i.e. $g(x) = \frac 1 { 2 \epsilon} {\bf 1}_ {(0,2 \epsilon )} (x).$
\par\noindent
To prove this, it suffices to verify that $q, \ \pi _0(x)$ given by \eqref{pi0forexampleuno}, and $g(x)$ satisfy Eq.  \eqref{maineq}, 
with $a=0$ and $b=2 \epsilon .$

\bigskip

\noindent {\bf Acknowledgments}\par\noindent
%The author would like to express particular thanks to the anonymous reviewer for his/her
%useful comments, leading to improved presentation. \par\noindent
The author acknowledges the MIUR Excellence Department Project awarded to
the Department of Mathematics, University of Rome Tor Vergata, CUP
E83C18000100006


\begin{thebibliography}{99}

\bibitem {abundo:saa20}
Abundo, M., 2020. \newblock  An inverse problem for the first-passage place of some diffusion processes with random starting point.
\newblock {Stochastic Anal. Appl.} vol 38, No 6, 1122--1133, https://doi.org/10.1080/07362994.2020.1768867


\bibitem {abundo:saa19}
Abundo, M., 2019. \newblock  An inverse first-passage problem revisited: the case of fractional Brownian motion, and time-changed Brownian motion.
\newblock {Stochastic Anal. Appl.}  vol 37, No 5, 708--716,  https://doi.org/10.1080/07362994.2019.1608834

\bibitem {abundo:mathematics}
Abundo, M., 2018. \newblock The Randomized First-Hitting Problem
of Continuously Time-Changed Brownian Motion.
\newblock {Mathematics} 6(6), 91, 1--10.
https://doi.org/10.3390/math6060091


\bibitem {abundo:LNSIM}
Abundo, M., 2015. \newblock An overview on inverse first-passage-time problems for one-dimensional diffusion processes.
\newblock {Lecture Notes of Seminario Interdisciplinare di Matematica} Vol. 12, 1 -- 44. http://dimie.unibas.it/site/home/info/documento3012448.html


\bibitem {abundo14b}
Abundo, M., 2014. \newblock One-dimensional reflected diffusions with two boundaries and an inverse first-hitting problem.
\newblock{Stochastic Anal. Appl.} 32, 975-991. DOI: 10.1080/07362994.2014.959595


\bibitem {abundo13b}
Abundo, M., 2013. \newblock Solving an inverse first-passage-time problem for Wiener process subject to
      random jumps from a boundary.
\newblock{Stochastic Anal. Appl.} 31: 4, 695-707.


\bibitem {abu:smj13}
Abundo, M., 2013. \newblock Some randomized first-passage problems for
one-dimensional diffusion processes.
\newblock{Scientiae Mathematicae Japonicae} 76(1), 33--46.

\bibitem  {abundo:stapro13}
Abundo, M., 2013. \newblock The double-barrier inverse first-passage problem for Wiener process with random
starting point.
\newblock{Stat. and Probab. Letters} 83, 168--176.

\bibitem  {abundo:stapro12}
Abundo, M., 2012. \newblock An inverse first-passage problem for one-dimensional diffusion with random starting point.
\newblock{Stat. and Probab. Letters} 82, 7--14. \ See also Erratum: Stat. and Probab. Letters, 82(3), 705.

\bibitem  {abundo:mcap10}
Abundo, M., 2010. \newblock On the First Hitting Time of a One-dimensional Diffusion and a Compound Poisson Process.
\newblock{Methodol. Comput. Appl.  Probab.} 12, 473--490.

\bibitem  {abundo:pms00}
Abundo, M., 2000. \newblock On first-passage-times for one-dimensional jump-diffusion processes.
\newblock{Prob. Math.Statis.} 20(2), 399--423.

\bibitem {gupta}
Gupta, A.K. (Ed.), Nadarajah, S. (Ed.), 2004. \newblock Handbook of Beta Distribution and Its Applications.
\newblock Boca Raton: CRC Press, https://doi.org/10.1201/9781482276596.

\bibitem {jackson:stapro09}
Jackson, K., Kreinin, A., and Zhang, W., 2009. \newblock
Randomization in the first hitting problem.
\newblock {Stat. and Probab. Letters} 79, 2422--2428.

%\bibitem {klebaner}
%Klebaner, F.C., 2005. \newblock Introduction to Stochastic Calculus with Applications, 2nd ed. \newblock
%London, Imperial College Press.


\bibitem {kouwang:03}
Kou, S.G. and Wang, H., 2003. \newblock First passage times of a jump diffusion process.
\newblock {Adv. Appl. Probab.} 35(2), 504--531.

\bibitem {lanska:89}
Lanska, V. and Smiths C.E., 1989. \newblock The effect of a random initial value in neural first-passage-time
models. \newblock {Math. Biosci.}  93, 191--215.

\bibitem {lefeb:20}
Lefebvre, M., 2020. \newblock First-passage problems for diffusion processes with state-dependent jumps.
\newblock {Communications in Statistics - Theory and Methods}, Published online: 02 Jul 2020.
https://doi.org/10.1080/03610926.2020.1784433


\bibitem {lefeb:19}
Lefebvre, M., 2019. \newblock Moments of First-Passage Places for Jump-Diffusion Processes.
\newblock {Sankhya A}, 1--9. https://doi.org/10.1007/s13171-019-00181-4

\bibitem {lefeb:19b}
Lefebvre, M., 2019. \newblock Minimizing the time spent in an interval by a Wiener process with uniform jumps.
\newblock {Control and Cybernetics} 48, No.3, 407--415.

\bibitem {tuckwell:76}
Tuckwell, H.C., 1976. \newblock On the first-exit time problem for temporally homogeneous Markov processes.
\newblock {J. Appl. Probab.} 13, 39--48.




\end{thebibliography}
\end{document}